%
%
%
\input amstex
\documentstyle{amsppt}

\NoBlackBoxes

\magnification=\magstep1
\hsize=6.5 true in
\vsize=8.2 true in
\topmatter
\title
Riemann-Stieltjes Integral Operators between Weighted Bergman Spaces
\endtitle

\author
                           {Jie Xiao}
\endauthor

\address
\fl Department of Mathematics and Statistics
\fl Memorial University of Newfoundland
\fl St. John's, NL A1C 5S7,  Canada
\fl Email: jxiao\@math.mun.ca
\endaddress

\thanks
\fl 2000 {\it Mathematics Subject Classification} 32H20, 30D55,
47B38 \fl This work was supported by NSERC of Canada. \fl Key words
and phrases: Sharpness; Boundedness; Compactness; Riemann-Stieltjes
integrals; Weighted Bergman spaces; $\Bbb C^n$-ball.
\endthanks

\def\qe{\gtrsim}

\def\fl{\flushpar}

\abstract
This note completely describes the bounded or compact Riemann-Stieltjes integral operators $T_g$ acting
between the weighted Bergman space pairs $(A^p_\alpha,A^q_\beta)$ in terms of particular regularities of
the holomorphic symbols $g$ on the open unit ball of $\Bbb C^n$.
\endabstract
\endtopmatter

\document
\leftheadtext{Jie Xiao}
\rightheadtext{Riemann-Stieltjes Integral Operators between Weighted Bergman Spaces}

\heading
{1. Introduction}
\endheading

Let ${\Bbb B}_n=\{z\in \Bbb C^n: \ |z|<1\}$ be the open unit ball of ${\Bbb C}^n$. Set $\partial\Bbb B_n=\{z\in\Bbb C^n:\ |z|=1\}$ be the compact unit sphere of ${\Bbb C}^n$ -- an $n$-dimensional Hilbert space over the complex field $\Bbb C$ under the inner product
$$
\langle z,w\rangle=\sum_{k=1}^n z_k\overline{w_k},\quad z=(z_1,z_2,...,z_n),\  w=(w_1,w_2,...,w_n)\in \Bbb C^n
$$
and the associated norm
$$
|z|=\langle z,z\rangle ^{1/2},\quad z\in\Bbb C^n.
$$
Given a holomorphic map $g:\Bbb B_n\to \Bbb C$. For a holomorphic map $f: \Bbb B_n\to \Bbb C$, we, as studied in [Hu], define the Riemann-Stieltjes integral of $f$ with respect to $g$ via
$$
T_g f(z)=\int_0^1 f(tz)\Cal R g(tz)t^{-1}{dt},\quad z\in \Bbb B_n\tag 1.1
$$
where
$$
{\Cal R} g(z)=\sum_{j=1}^n z_j\frac{\partial g(z)}{\partial z_j},\quad z=(z_1,z_2,...,z_n)\in\Bbb B_n
$$
stands for the radial derivative of $g$. In particular, if
$$
{\bold 1}=(1,1,...,1)\quad\hbox{and}\quad g(z)=-\log(1-\langle z,{\bold 1}\rangle)$$
then (1.1) becomes
$$
T_g f(z)=\int_0^1 f(tz)\langle tz,{\bold 1}\rangle (1-\langle tz,{\bold 1}\rangle)^{-1}t^{-1}{dt},\quad z\in \Bbb B_n,
$$
a higher dimensional version of the classical Ces\'aro operator.

We consider the problem of determining the optimum on $g$ such that
$$
T_g: {A^p_\alpha}\to {A^q_\beta}\quad {boundedly\quad or\quad compactly}.\tag 1.2
$$
Here and henceforth, for $p>0$ and $\alpha>-1$, $A^p_\alpha$ is the weighted Bergman space of all holomorphic maps $f:\Bbb B_n\to \Bbb C$ satisfying
$$
\|f\|_{A^p_\alpha}=\left(\int_{\Bbb B_n}|f(z)|^p(1-|z|^2)^\alpha dv(z)\right)^{1/p}<\infty,
$$
where $dv$ denotes the Lebesgue volume measure on $\Bbb B_n$.
It is known that if $p=q\ge 1$ and $\alpha=\beta>-1$ then (1.2) holds if and only if $g$ belongs to the Bloch space or the little Bloch space; see for example [AlSi] (for $n=1$) and [Xi] (for $n\ge 1$). But, in other cases (even for $n=1$), the optimal $(A_\alpha^p, A^q_\beta)$ estimates have not yet been worked out; see also [AlCi] (and references therein) for the setting of $(H^p, H^q)$--the $1$-dimensional limit case of $(A^p_\alpha, A^q_\beta)$ as $\alpha=\beta\to -1$.

The goal of this note is to find out the optimal conditions of $g$ such that (1.2) holds in all possible cases on $(p,q)\in (0,\infty)\times(0,\infty)$ and $(\alpha,\beta)\in (-1,\infty)\times(-1,\infty)$. Below is our result.

\proclaim{Theorem 1.1} Let $\alpha,\beta\in (-1,\infty)$, $p,q\in (0,\infty)$ and $g:\Bbb B_n\to \Bbb C$ be holomorphic.

\flushpar{\rm(i)} If $p>q$ then $T_g: A^p_\alpha\to A^q_\beta$ boundedly or compactly if and only if
$$
\int_{\Bbb B_n}\big(|\Cal R g(z)|(1-|z|^2)^{1-\frac{\alpha}{p}+\frac{\beta}{q}}\big)^{\frac{pq}{p-q}} dv(z)<\infty. \tag 1.3
$$
\flushpar{\rm(ii)} If $p\le q$ then $T_g: A^p_\alpha\to A^q_\beta$ boundedly or compactly if and only if
$$
\sup_{z\in {\Bbb B_n}}|\Cal R g(z)| (1-|z|^2)^{1-\frac{n+1+\alpha}{p}+\frac{n+1+\beta}{q}}<\infty\tag 1.4
$$
or
$$
\lim_{z\to\partial{\Bbb B_n}}|\Cal R g(z)| (1-|z|^2)^{1-\frac{n+1+\alpha}{p}+\frac{n+1+\beta}{q}}=0. \tag 1.5
$$
\endproclaim

\medskip

It is worth pointing out that the symbol $g$ satisfying (1.3) or (1.4)/(1.5) is a constant whenever
$$
1\le\frac{1+\alpha}{p}-\frac{1+\beta}{q}\ \hbox{and}\ p>q\quad\hbox{or}\quad \frac{n+1+\alpha}{p}-\frac{n+1+\beta}{q}>1\ \hbox{and}\ p\le q.
$$

The rest of this note is organized as follows. In Section 2, we collect some preliminary but useful facts on the weighted Bergman spaces and Khinchine's inequality. In Section 3, we demonstrate Theorem 1.1 through these preliminary results and some of the ideas exposed in [Lu] and [SmYa].

\heading
{2. Preliminaries}
\endheading

When $p\ge 1$, the space $A^p_\alpha$ is a Banach space equipped with the norm $\|\cdot\|_{A^p_\alpha}$, and when $p\in (0,1)$, the space $A^p_\alpha$ is a complete metric space with the distance $d(f,g)=\|f-g\|^p_{A^p_\alpha}.$

First of all, we need a growth property of holomorphic functions on $\Bbb B_n$. To do so, we denote by $\phi_w$ the automorphism of $\Bbb B_n$ taking $0$ to $w\in\Bbb B_n$, i.e.,
$$
\phi_w(z)=\cases z, & w={0}\\
\frac{w-P_wz-\sqrt{1-|z|^2}Q_wz}{1-\langle z,w\rangle}, & w\not={0},
\endcases
$$
where $Q_w=I-P_w$, $I$ is the identity map and $P_w$ is the projection of $\Bbb C^n$ onto the one-dimensional subspace spanned by $w\not=0$. For $r>0$ and $z\in\Bbb B_n$ put
$$
D(z,r)=\left\{w\in\Bbb B_n: \frac{1}{2}\log\frac{1+|\phi_w(z)|}{1-|\phi_w(z)|}<r\right\}
$$
which is called the Bergman metric ball with center $z$ and radius $r$.

\proclaim{Lemma 2.1} Let $p\in (0,\infty)$ and $\alpha\in (-1,\infty)$. If $f:\Bbb B_n\to \Bbb C$ is holomorphic then
$$
|f(z)|(1-|z|^2)^{\frac{n+1+\alpha}{p}}\lesssim \left(\int_{D(z,r)}|f(w)|^p(1-|w|^2)^\alpha dv(w)\right)^{\frac{1}{p}},\quad z\in\Bbb B_n.\tag 2.1
$$
\endproclaim
For a proof of (2.1), see also [Zh, p. 64].

Next, we state two characterizations of the weighted Bergman spaces.

\proclaim{Lemma 2.2} Let $p\in (0,\infty)$, $\alpha\in (-1,\infty)$, and $f:\Bbb B_n\to \Bbb C$ be holomorphic. Then the following three statements are equivalent:
\smallskip
\flushpar{\rm(i)} $f\in A^p_\alpha$.
\smallskip
\flushpar{\rm(ii)}
$$
|\!|\!|f|\!|\!|_{A^p_\alpha}=|f({0})|+\left(\int_{\Bbb B_n}|\Cal R f(z)|^p(1-|z|^2)^{p+\alpha}dv(z)\right)^{\frac{1}{p}}<\infty.
$$
\flushpar{\rm (iii)} For any $\eta\in (0,1]$, there exist a sequence $\{z_j\}$ in $\Bbb B_n$  such that
\smallskip
{\rm (a)} $\Bbb B_n=\cup_j D(z_j,\eta)$;
\smallskip
{\rm (b)} $D(z_j,\frac{\eta}{4})\cap D(z_k,\frac{\eta}{4})=\emptyset$ for $j\not=k$;
\smallskip
{\rm (c)} Each point $z\in \Bbb B_n$ lies in at most $N=N(\eta)$ of balls from $\{D(z_j, 2\eta)\}$;
\smallskip
{\rm (d)}
$$
f(z)=\sum_j c_j \frac{(1-|z_j|^2)^{\frac{pb-n-1-\alpha}{p}}}{(1-\langle z,z_j\rangle)^b},\quad z\in\Bbb B_n,
$$
where $\{c_j\}$ is in the sequence space $l^p$ and $b$ is a constant greater than $n\max\{1,\frac{1}{p}\}+\frac{1+\alpha}{p}.$
\smallskip
\flushpar Moreover, if $f\in A^p_\alpha$ then
\smallskip
$$
|\!|\!|f|\!|\!|_{A^p_\alpha}\approx\|f\|_{A^p_\alpha}\approx\|\{c_j\}\|_{l^p}.\tag 2.2
$$
\endproclaim

For a proof of Lemma 2.2 and its sources, see for example [Zh, Chapter 2] and [CoRo] as well as [Sh].

Finally, we quote the following well-known form of Khinchine's inequality.

\proclaim{Lemma 2.3} Suppose
$$
r_0(t)=\cases 1, & 0\le t-[t]<1/2\\
-1, & 1/2\le t-[t]<1
\endcases
$$
and
$$
r_j(t)=r_0(2^jt), \quad j=1,2,...,
$$
and let $p\in (0,\infty)$ and $(c_1,...,c_m)\in \Bbb C^m$, $m=1,2,....$ Then
$$
\left(\sum_{j=1}^m |c_j|^2\right)^{\frac{1}{2}}\approx\left(\int_0^1\Big|\sum_{j=1}^m c_j r_j(t)\Big|^p dt\right)^{\frac{1}{p}}.\tag 2.3
$$
\endproclaim

{\it Note}: In the above (and below), the notation $U\approx V$ means that there are two constants $\kappa_1, \kappa_2>0$ such that $\kappa_1 V\le U\le \kappa_2 V$. Moreover, if $U\le \kappa_2 V$ then we say $U\lesssim V$.

\heading
{3. Proof}
\endheading

To begin with, we notice that the formula
$$
\Cal R T_g f(z)=f(z) \Cal R g(z),\quad z\in\Bbb B_n,\tag 3.1
$$
holds for all holomorphic maps $f, g:\Bbb B_n\to \Bbb C$. And, let us agree to two more conventions:
$$
dv_\beta(z)=(1-|z|^2)^\beta dv(z),\quad z\in \Bbb B_n,
$$
and
$$
\|T_g\|_{A^p_\alpha\to A^q_\beta}=\sup_{\|f\|_{A^p_\alpha}=1}\|T_g f\|_{A^q_\beta}.
$$
It is clear that
$$
\|T_g f\|_{A^q_\beta}\le\|T_g\|_{A^p_\alpha\to A^q_\beta}\|f\|_{A^p_\alpha},\quad f\in A^p_\alpha.\tag 3.2
$$

\medskip

\demo{\bf Proof of Theorem 1.1 (i)} Suppose $p>q$. Let (1.3) be true. By (2.2), (3.1) and H\"older's inequality we have that if $f\in A^p_\alpha$ then

$$
\align
\|T_g f\|^q_{A^q_\beta}&\approx\int_{\Bbb B_n}|f(z)\Cal R g (z)|^q(1-|z|^2)^q dv_\beta(z)\\
&\lesssim\|f\|^q_{A^p_\alpha}\left(\int_{\Bbb B_n}\big(|\Cal R g(z)|(1-|z|^2)^{1-\frac{\alpha}{p}+\frac{\beta}{q}}\big)^{\frac{pq}{p-q}} dv(z)\right)^{\frac{p-q}{p}},
\endalign
$$
implying the boundedness of $T_g: A^p_\alpha\to A^q_\beta$.

Conversely, suppose $T_g: A^p_\alpha\to A^q_\beta$ is bounded. Then $\|T_g\|_{A^p_\alpha\to A^q_\beta}$ is finite with (3.2). For each natural number $j$ let
$$
K_j(z)=\frac{(1-|z_j|^2)^{\frac{pb-n-1-\alpha}{p}}}{(1-\langle z,z_j\rangle)^b},
$$
where $\{z_j\}$ is the sequence in Lemma 2.2 (iii). Let $\{c_j\}\in l^p$, and choose $\{r_j(t)\}$ as obeying Lemma 2.3. Then $\{c_jr_j(t)\}\in l^p$ with
$\|\{c_jr_j(t)\}\|_{l^p}=\|\{c_j\}\|_{l^p},$ and so $\sum_j c_j r_j(t) K_j\in A^p_\alpha$ with
$$
\Big\|\sum_j c_j r_j(t) K_j\Big\|_{A^p_\alpha}\approx\|\{c_j\}\|_{l^p},
$$
due to Lemma 2.2 (iii). This fact plus (3.2) derives
$$
\int_{\Bbb B_n}\Big|T_g\Big(\sum_j c_j r_j(t)K_j\Big)(z)\Big|^q dv_\beta(z)\lesssim\|T_g\|_{A^p_\alpha\to A^q_\beta}^q\|\{c_j\}\|_{l^p}^q,
$$
Furthermore, integrating this inequality from $0$ to $1$ with respect to $dt$, as well as using (3.1), Fubini's theorem and (2.3) in Lemma 2.3, we get
$$
\int_{\Bbb B_n}\Big(\sum_j \big|c_j K_j(z)\big|^2\Big)^{\frac{q}{2}}\big(|\Cal R g(z)|(1-|z|^2)\big)^q dv_\beta(z)\lesssim\|T_g\|_{A^p_\alpha\to A^q_\beta}^q\|\{c_j\}\|_{l^p}^q.
$$
Noticing the estimate
$$
|1-\langle z,z_j\rangle|\approx 1-|z_j|^2\quad\hbox{as}\quad z\in D(z_j,2\eta),
$$
applying the condition (c) in Lemma 2.2 (iii), and letting $1_E$ be the characteristic function of a set $E\subseteq \Bbb B_n$, we achieve
$$
\align
\sum_j&|c_j|^q\left(\frac{\int_{D(z_j,2\eta)}\big(|\Cal R g(z)|(1-|z|^2)\big)^q dv_\beta(z)}{(1-|z_j|^2)^{\frac{q(n+1+\alpha)}{p}}}\right)\\
&=\int_{\Bbb B_n}\sum_j\left(\frac{|c_j|^q1_{D(z_j,2\eta)}(z)}{
(1-|z_j|^2)^{\frac{q(n+1+\alpha)}{p}}}\right)\big(|\Cal R g(z)|(1-|z|^2)\big)^q dv_\beta(z)\\
&\lesssim\int_{\Bbb B_n}\left(\sum_j\frac{|c_j|^21_{D(z_j,2\eta)}(z)}{(1-|z_j|^2)^{\frac{2(n+1+\alpha)}{p}}}\right)^{\frac{q}{2}}\big(|\Cal R g(z)|(1-|z|^2)\big)^q dv_\beta(z)\\
&\lesssim \int_{\Bbb B_n}\Big(\sum_j|c_jK_j(z)|^2\Big)^{\frac{q}{2}}\big(|\Cal R g(z)|(1-|z|^2)\big)^q dv_\beta(z)\\
&\lesssim\|T_g\|_{A^p_\alpha\to A^q_\beta}^q\|\{c_j\}\|_{l^p}^q.
\endalign
$$
The last estimate indicates
$$
\left\{\frac{\int_{D(z_j,2\eta)}\big(|\Cal R g(z)|(1-|z|^2)\big)^q dv_\beta(z)}{(1-|z_j|^2)^{\frac{q(n+1+\alpha)}{p}}}\right\}\in l^{\frac{p}{p-q}}.
$$
Because $\Cal R g$ is holomorphic on $\Bbb B_n$, by (2.2), (2.1) and the condition (a) in Lemma 2.2 (iii) we achieve
$$
\align
\int_{\Bbb B_n}&\big(|\Cal R g(z)|(1-|z|^2)\big)^{\frac{pq}{p-q}} (1-|z|^2)^{\frac{\beta p-\alpha q}{p-q}}dv(z)\\
&\lesssim\sum_j\int_{D(z_j,\eta)}\big(|\Cal R g(z)|(1-|z|^2)\big)^{\frac{pq}{p-q}} (1-|z|^2)^{\frac{\beta p-\alpha q}{p-q}}dv(z)\\
&\lesssim\sum_j\int_{D(z_j,\eta)}\left(\frac{\int_{D(z,\eta)}\big(|\Cal R g(w)|(1-|w|^2)\big)^q dv_\beta(w)}{(1-|z|^2)^{\frac{q\alpha+p(n+1)}{p}}}\right)^{\frac{p}{p-q}}dv(z)\\
&\lesssim\sum_j\left(\frac{\int_{D(z_j,2\eta)}\big(|\Cal R g(z)|(1-|z|^2)\big)^q dv_\beta(z)}{(1-|z_j|^2)^{\frac{q\alpha+p(n+1)}{p}}}\right)^{\frac{p}{p-q}}(1-|z_j|^2)^{n+1}\\
&\lesssim\sum_j\left(\frac{\int_{D(z_j,2\eta)}\big(|\Cal R g(z)|(1-|z|^2)\big)^q dv_\beta(z)}{(1-|z_j|^2)^{\frac{q(n+1+\alpha)}{p}}}\right)^{\frac{p}{p-q}}\\
&\lesssim \|T_g\|^{\frac{pq}{p-q}}_{A^p_\alpha\to A^q_\beta},
\endalign
$$
giving (1.3).

Regarding the compactness, it suffices to show that if (1.3) holds then $T_g: A^p_\alpha\to A^q_\beta$ is compact. Assuming (1.3), we obtain that $T_g$ is a bounded operator from $A_\alpha^p\to A_\beta^q$ and so $g=g({0})+T_g 1\in A^q_\beta$. In addition to this, we also see that for any $\epsilon>0$ there is a $\delta\in (0,1)$ such that
$$
\int_{|z|>\delta}\Big(|\Cal R g(z)|(1-|z|^2)\Big)^{\frac{pq}{p-q}}(1-|z|^2)^{\frac{\beta p-\alpha q}{p-q}}dv(z)<\epsilon.
$$
Since the weak convergence in $A^p_\alpha$ means the uniform convergence on compacta of $\Bbb B_n$, we may assume that $\{f_j\}$ is any sequence in the unit ball of $A_\alpha^p$ and converges to $0$ uniformly on compacta of $\Bbb B_n$. For the above $\epsilon>0$ there exists an integer $j_0>0$ such that $\sup_{|z|\le\delta}|f_j(z)|<\epsilon$ as $j\ge j_0$. With the help of (3.1), (2.2) and H\"older's inequality, we further obtain
$$
\align
\|T_g f_j\|^q_{A_\beta^q}&\approx\left(\int_{|z|\le\delta}+\int_{|z|>\delta}\right)\Big(|f_j(z)\Cal R g(z)|(1-|z|^2)\Big)^q dv_\beta(z)\\
&\lesssim \epsilon^q \|g\|_{A^q_\beta}^q+\|f_j\|_{A^p_\alpha}^q\left(\int_{|z|>\delta}\Big(|\Cal R g(z)|(1-|z|^2)\Big)^{\frac{pq}{p-q}}(1-|z|^2)^{\frac{\beta p-\alpha q}{p-q}}dv(z)\right)^{\frac{p-q}{p}}\\
&\lesssim \epsilon^q\|g\|_{A^q_\beta}^q+\epsilon^{\frac{p-q}{p}}.
\endalign
$$
In other words, $\lim_{j\to\infty}\|T_g f_j\|_{A_\beta^q}=0$ and so $T_g: A^p_\alpha\to A^q_\beta$ is compact.
\enddemo

\medskip
\demo{\bf Proof of Theorem 1.1 (ii)} Suppose now $p\le q$. If (1.4) holds then
\medskip

$$
\|g\|_{B^\gamma}=\sup_{z\in {\Bbb B_n}}|\Cal R g(z)|(1-|z|^2)^{\gamma}<\infty,\quad \gamma=1-\frac{n+1+\alpha}{p}+\frac{n+1+\beta}{q}.
$$
From (3.1), (2.2) and (2.1) it turns out that for $f\in A^p_\alpha$,
$$
\align
\|T_g f\|_{A^q_\beta}^q&\approx\int_{\Bbb B_n}|f(z)|^p|f(z)|^{q-p}|\Cal R g(z)|^q(1-|z|^2)^{q}dv_\beta(z)\\
&\lesssim\|f\|_{A^p_\alpha}^{q-p}\|g\|_{B^\gamma}^q\int_{\Bbb B_n}|f(z)|^p(1-|z|^2)^{q-q\gamma-\frac{(q-p)(n+1+\alpha)}{p}}dv_\beta(z)\\
&\lesssim\|f\|_{A^p_\alpha}^{q}\|g\|_{B^\gamma}^q.
\endalign
$$
That is to say, $T_g: A^p_\alpha\to A^q_\beta$ is bounded.

Conversely, if $T_g: A^p_\alpha\to A^q_\beta$ is bounded then the operator norm $\|T_g\|_{A^p_\alpha\to A^q_\beta}$ is finite with (3.2). Keeping this in mind, we deal with two cases: $p>1$ and $p\le 1$.
\medskip

{\it Case 1}: $p>1$. Define
$$
K(w,z)={(1-\langle z,w\rangle)^{-(n+\alpha+1)}},\quad z,w\in\Bbb B_n.
$$
A routine calculation (see for example [FaKo] or [Zh, pp. 20-21]) yields
$$
\|K(w,\cdot)\|_{A^p_\alpha}\approx{(1-|w|^2)^{-\frac{(n+1+\alpha)(p-1)}{p}}}.\tag 3.3
$$
Using (3.1), (2.2), certain transformation properties of $\phi_w$ and (2.1) (for $\Cal R g(\phi_w)$), we obtain
$$
\align
\|T_g &K(w,\cdot)\|_{A^q_\beta}^q\approx\int_{\Bbb B_n}|K(w,z)|^q|\Cal R g(z)|^q(1-|z|^2)^{q}dv_\beta(z)\\
&\qe (1-|w|^2)^{q+(1-q)(n+1+\alpha)+\beta-\alpha}\int_{|u|\le 1/2}
\frac{|\Cal R g(\phi_w(u))|^q(1-|u|^2)^{q+\beta}}{|1-\langle u,w\rangle|^{2n+2-(n+1+\alpha)q+2(\beta-\alpha)}}dv(u)\\
&\qe (1-|w|^2)^{q+(1-q)(n+1+\alpha)+\beta-\alpha}|\Cal R g(\phi_w(0))|^q\\
&\qe (1-|w|^2)^{q+(1-q)(n+1+\alpha)+\beta-\alpha}|\Cal R g(w)|^q.
\endalign
$$
In brief, we have
$$
\|T_g K(w,\cdot)\|_{A^q_\beta}^q\qe (1-|w|^2)^{q+(1-q)(n+1+\alpha)+\beta-\alpha}|\Cal R g(w)|^q.\tag 3.4
$$
This estimate, together with (3.3) and (3.2) (for $f(\cdot)=K(w,\cdot)$), produces (1.4) right away.

{\it Case 2}: $p\le 1$. Select a positive integer $m>{n+1+\alpha}$ and set
$$
K_{p}(w,z)={\big(1-\langle z,w\rangle\big)^{-\frac{m}{p}}},\quad z,w\in\Bbb B_n.
$$
Just like the case of $p>1$, it follows that
$$
\|K_{p}(w,\cdot)\|_{A^p_\alpha}
\approx (1-|w|^2)^{\frac{n+1+\alpha-m}{p}}\tag 3.5
$$
and
$$
\|T_g K_{p}(w,\cdot)\|^q_{A^q_\beta}
\qe(1-|w|^2)^{q+n+1+\beta-\frac{mq}{p}}|\Cal R g(w)|^q.\tag 3.6
$$
A combination of (3.6), (3.5) and (3.2) (for $f(\cdot)=K_p(w,\cdot)$) yields (1.4) too.

To establish the corresponding compactness part, we assume that $g$ satisfies (1.5). Then $g\in A_\beta^q$, and for any $\epsilon>0$ there is an $\delta\in (0,1)$ such that as $|z|\in (\delta,1)$,
$$
|\Cal R g(z)|(1-|z|^2)^{\gamma}<\epsilon,\quad \gamma=1-\frac{n+1+\alpha}{p}+\frac{n+1+\beta}{q}.
$$
In order to prove that $T_g: A^p_\alpha\to A^q_\beta$ is compact, we consider any sequence $\{f_j\}$ in the unit ball of $A^p_{\alpha}$ which converges to $0$ uniformly on compacta of $\Bbb B_n$. For such a sequence, there is an integer $j_0>0$ such that $\sup_{|z|\le \delta}|f_j(z)|<\epsilon$ when $j\ge j_0.$
Hence by (3.1), (2.1) and (2.2),
$$
\align
\|T_g f_j\|^q_{A^q_\beta}
&\approx\int_{\Bbb B_n}|f_j(z)\Cal R g(z)|^q(1-|z|^2)^{q}dv_\beta(z)\\
&\lesssim\epsilon^q\int_{|z|\le \delta}|\Cal R g(z)|^q(1-|z|^2)^{q}dv_\beta(z)+\epsilon^q\|f_j\|_{A^p_\alpha}^{q-p}\int_{|z|>\delta}|f_j(z)|^p (1-|z|^2)^\alpha dv(z)\\
&\lesssim\epsilon^q\Big(\|g\|_{A^q_\beta}^q+\|f_k\|_{A^p_\alpha}^q\Big)\\
&\lesssim\epsilon^q\Big(\|g\|_{A^q_\beta}^q+1\Big).
\endalign
$$
Namely, $\|T_g f_j\|_{A^q_\beta}\to 0$ {as} $j\to\infty.$ Therefore $T_g$ is a compact operator from $A^p_\alpha$ to $A^q_\beta$.

On the other hand, if $T_g: A^p_\alpha\to A^q_\beta$ is compact, then, for $z,w\in\Bbb B_n$ let
$$
k_p(w,z)=\cases \frac{K(w,z)}{\|K(w,\cdot)\|_{A^p_\alpha}},\quad & p>1\\
\frac{K_{p}(w,z)}{\|K_{p}(w,\cdot)\|_{A^p_\alpha}},\quad & p\le 1.
\endcases
$$
Obviously, $k_p(w,\cdot)$ tend to $0$ uniformly on compacta of $\Bbb B_n$ as $w\to\partial\Bbb B_n$. By the compactness of $T_g: A^p_\alpha\to A^q_\beta$, we find
$$
\lim_{w\to\partial\Bbb B_n}\|T_g k_p(w,\cdot)\|_{A^q_\beta}=0.
$$
The above limit, along with (3.3), (3.4), (3.5) and (3.6), yield (1.5).
\enddemo

\medskip
\flushpar{\it Acknowledgment}. The author is grateful to K. Zhu for his helpful emails.
\medskip

\heading
 {References}
\endheading

\flushpar {[AlCi]}\ A. Aleman and J. A. Cima, {An integral operator on $H^p$ and Hardy's inequality}, {\it J. Anal. Math.}, {85}(2001), 157-176.
\smallskip

\flushpar {[AlSi]}\ A. Aleman and A. G. Siskakis, {Integral operators on Bergman spaces}, {\it Indiana Univ. Math. J.}, {46}(1997), 337-356.
\smallskip

\flushpar {[CoRo]}\ R. Coifman and R. Rochberg, {Representation theorems for holomorphic and harmonic functions in $L^p$}, {\it Asterisque}, 77(1980), 157-165.
\smallskip

\flushpar {[FaKo]}\ J. Faraut and A. Koranyi, {Function spaces and reproducing kernels on bounded symmetric domains}, {\it J. Funct. Anal.}, {88} (1990), 64-89.
\smallskip

\flushpar {[Hu]}\ Z. Hu, {Extended Ces\'aro operators on the Bloch space in the unit ball of $\Bbb C^n$}, {\it Acta Math. Sci. Ser. B Engl. Ed.}, {23}(2003), 561-566.
\smallskip

\flushpar {[Lu]}\ D. H. Luecking, {Embedding theorems for spaces of analytic functions via Khinchine's inequality}, {\it Michigan Math. J.}, {40} (1993), 333-358.
\smallskip

\flushpar {[Sh]}\ J. Shapiro, {Macey topologies, reproducing kernels, and diagonal maps on the Hardy and Bergman spaces}, {\it Duke Math. J.}, {43}(1976), 187-202.
\smallskip

\flushpar {[SmYa]}\ W. Smith and L. Yang, {Composition operators that improve integrability on weighted Bergman spaces}, {\it Proc. Amer. Math. Soc.}, {126} (1998), 411-420.
\smallskip

\flushpar {[Xi]}\ J. Xiao, {Riemann-Stieltjes operators on weighted Bloch and Bergman spaces of the unit ball}, {\it J. London. Math. Soc.}, (2) {70} (2004), 199-214.
\smallskip

\flushpar {[Zh]}\ K. Zhu, {\it Spaces of Holomorphic Functions in the Unit Ball}, Graduate Texts in Mathematics, 226. Springer-Verlag, New York, 2005.
\smallskip

\smallskip

\enddocument